\documentclass[a4paper,10pt,
]{article}
\usepackage{amsmath,amssymb}
\usepackage[dvips]{color}
\usepackage{cite}
\usepackage{hangcaption}
\usepackage{amsmath,amssymb}
\usepackage[dvips]{color}
\usepackage{cite}
\usepackage{eepic}
\usepackage{epic}
\usepackage{hangcaption}
\usepackage{amsmath}
\usepackage{amssymb}
\usepackage{amscd}
\usepackage{amsthm}
\usepackage[dvipdfmx,%
bookmarks=true,%
bookmarksnumbered=true,%
setpagesize=false,%
pdfkeywords={TeX; dvipdfmx; hyperref; color;}]{hyperref}
\usepackage{hyperref}
\def\cases{\left\{\begin{array}{ll}}
\def\endcases{\end{array}\right.}
\def\roman{\rm}
\def\bigtimes{\mathop{\mbox{\Large $\times$}}}
\begin{document}
\setcounter{page}{1}
\vskip1.5cm
\begin{center}
{\LARGE \bf 
Kalman filter in quantum language}
\vskip0.5cm
{\rm
\large
Shiro Ishikawa,  Kohshi Kikuchi
}
\\
\vskip0.2cm
\rm
\it
\it
Department of Mathematics, Faculty of Science and Technology,
Keio University,
\\ 
3-14-1, Hiyoshi, Kouhoku-ku Yokohama, Japan.
{E-mail:
ishikawa@math.keio.ac.jp,}%
\newline
\newline
Department of Applied Molecular Bioscience, Graduate School of Medicine, Yamaguchi University
\\
2-16-1 Tokiwadai, Ube 755-8611, Japan,
E-mail:
kohshi.kikuchi@gmail.com
\end{center}
\par
\rm
\vskip0.4cm
\par
\noindent
{\bf Abstract}
\normalsize
\par
\noindent
Recently, we proposed measurement theory ( or. quantum language) as a  linguistic turn of quantum mechanics (with the Copenhagen interpretation). This theory has a great power of scientific descriptions. In fact, we have continued asserting that even statistics can be described in terms of measurement theory. Thus, we believe that quantum language is future statistics (i.e., statistics will develop into quantum language). However, now we think that our arguments were too abstract and philosophical, that is, we should have presented concrete examples much more. Thus, in this paper, we show that the calculation of Kalman filter is more understandable in terms of quantum language than in terms of usual statistics. For this, we devote ourselves to statistical measurement theory, in which the Bertrand paradox is discussed.

\par
\vskip1.0cm
\par

\par
\def\Cal{\cal}
\def\bigstimes{\text{\large $\: \boxtimes \,$}}

\section{Measurement Theory
(Axioms
and
Interpretation)}

\subsection{The classifications of measurement theory}

\rm
\par
\par
\noindent
\par
In this section,
we introduce measurement theory (or in short, MT). This theory is a kind of language,
and thus, it is also called quantum language
(or in short, QL).

\par
\par
\rm
Measurement theory
({\it cf.} refs.\textcolor{black}{\cite{Ishi3}-\cite{Kiku1}})
is,
by an analogy of
quantum mechanics
(or,
as a linguistic turn of quantum
mechanics
), constructed
as the scientific
theory
formulated
in a certain 
{}{$C^*$}-algebra ${\cal A}$
(i.e.,
a norm closed subalgebra
in the operator algebra $B(H)$
composed of all bounded linear operators on a Hilbert space $H$,
{\it cf.$\;$}\textcolor{black}{\cite{Neum, Saka}}
). Let ${\mathcal N}$ be the weak${}^\ast$ closure of ${\mathcal A}$,
which is called a $W^\ast$-algebra. The structure
$[{\mathcal A} \subseteq {\mathcal N} \subseteq B(H)]$
is called a fundamental structure of MT.
\par
\noindent
\par
MT (= measurement theory) is composed of
two theories
(i.e.,
pure measurement theory
(or, in short, PMT]
and
statistical measurement theory
(or, in short, SMT).
That is,
we have:
\par
\rm
\par
\begin{itemize}
\item[(A$_1$)]
$
\underset{\text{\footnotesize }}{
\text{
MT 
}
}
\cases
\text{[PMT](related to Fisher statistics)}
\underset{\text{\scriptsize }}{\text{
}}
=
\displaystyle{
{
\mathop{\mbox{[(pure) measurement]}}_{\text{\scriptsize (Axiom$^{\rm P}$ 1) }}
}
}
+
\displaystyle{
\mathop{
\mbox{
[causality]
}
}_{
{
\mbox{
\scriptsize
(Axiom 2)
}
}
}
}
\\
\\
\text{[SMT](related to Baysian statitics)}
\underset{\text{\scriptsize }}{\text{
}}
=
\displaystyle{
{
\mathop{\mbox{[(statistical) measurement]}}_{\text
{\scriptsize (Axiom$^{\roman S}$ 1) }}
}
}
\!
+
\!
\displaystyle{
\mathop{
\mbox{
[causality]
}
}_{
{
\mbox{
\scriptsize
(Axiom 2)
}
}
}
}
\endcases
$
\end{itemize}
where
Axiom 2 is common in PMT and SMT.
For completeness, note that measurement theory (A$_1$)
is
a kind of language
based on
{\lq\lq}the quantum mechanical world view{\rq\rq},
({\it cf.$\;$} refs. \textcolor{black}{\cite{Ishi3, Ishi7, Ishi8}}).
%

When ${\cal A}=B_c(H)$,
the ${C^*}$-algebra composed
of all compact operators on a Hilbert space $H$,
the MT is called {quantum measurement theory}
(or,
quantum system theory),
which can be regarded as
the linguistic aspect of quantum mechanics.
Also, when ${\cal A}$ is commutative
$\big($
that is, 
when ${\cal A}$ is characterized by $C_0(\Omega)$,
the $C^*$-algebra composed of all continuous 
complex-valued functions vanishing at infinity
on a locally compact Hausdorff space $\Omega$
({\it cf.$\;$}\textcolor{black}{\cite{Saka}})$\big)$,
the MT is called {classical measurement theory}.
Thus, we have the following classification:
\begin{itemize}
\item[(A$_2$)]
$
\quad
\underset{\text{\scriptsize }}{\text{MT}}
$
$\left\{\begin{array}{ll}
\text{quantum MT$\quad$(when non-commutative ${\cal A}=B_c (H)$,
${\mathcal N}=B(H)$)}
%
\\
\\
\text{classical MT
$\quad$
(when commutative ${\cal A}=C_0(\Omega)$),
${\mathcal N}=L^\infty (\Omega, \nu)(\subseteq B(L^2(\Omega, \nu) )$)}
\end{array}\right.
$
\end{itemize}
Also,
MT has two formulations as follows.
\begin{itemize}
\item[(A$_3$)]
$
\quad
\underset{\text{\scriptsize }}{\text{MT}}
$
$\left\{\begin{array}{ll}
\text{MT$^{C^*}$($C^*$-algebraic formulation); the $C^*$-algebra ${\mathcal A}$ plays an important role
}
%
\\
\\
\text{MT$^{W^*}$($W^*$-algebraic formulation); the $W^*$-algebra ${\mathcal N}$ plays an important role
}
\end{array}\right.
$
\end{itemize}
In this paper, we devote ourselves to classical SMT${}^{W^*}$
(i.e.,
the classical statistical measurement theory with the $W^*$-algebra formulation
).

\subsection{Observables}
\par
\noindent
\par
Now we shall explain the measurement theory
(i.e.,
classical SMT${}^{W^*}$
).
\par
\noindent
\par
Let
$[{\mathcal A} \subseteq {\mathcal N} \subseteq B(H)]$
be the fundamental structure
of measurement theory.
Let
${\cal N}_*$ be the
pre-dual Banach space of
${\cal N}$.
That is,
$ {\cal N}_* $
$ {=}  $
$ \{ \rho \; | \; \rho$
is a weak$^*$ continuous linear functional on ${\cal N}$
$\}$,
and
the norm $\| \rho \|_{ {\cal N}_* } $
is defined by
$ \sup \{ | \rho ({}F{}) |  \:{}: \; F \in {\cal N}
\text{ such that }\| F \|_{{\mathcal N}} 
(=\| F \|_{B(H)} )\le 1 \}$.
The bi-linear functional
$\rho(F)$
is
also denoted by
${}_{{\cal N}_*}
\langle \rho, F \rangle_{\cal N}$,
or in short
$
\langle \rho, F \rangle$.
Define the
\it
mixed state
$\rho \;(\in{\cal N}_*)$
\rm
such that
$\| \rho \|_{{\mathcal N}_* } =1$
and
$
\rho ({}F) \ge 0
\text{ 
for all }F\in {\cal N}
\text{ satisfying }
F \ge 0$.
And put
\begin{align*} {\frak S}^m  ({}{\cal N}_*{})
{=}
\{ \rho \in {\cal N}_*  \; | \;
\rho
\text{ is a mixed state}
\}.
\end{align*}
%
\rm

According to the noted idea ({\it cf.} ref. \textcolor{black}{\cite{ Davi}})
in quantum mechanics,
an {\it observable}
${\mathsf O}{\; \equiv}(X, {\cal F},$
$F)$ in the $W^*$-algebra
${{\cal N}}$
is defined as follows:
\par
\par
\begin{itemize}
\item[(B$_1$)]
[$\sigma$-field]
$X$ is a set,
${\cal F}
(\subseteq 2^X $,
the power set of $X$)
is a $\sigma$-field of $X$,
that is,
{\lq\lq}$\Xi_1, \Xi_2, \Xi_3, \cdots \in {\cal F}\Rightarrow \cup_{k=1}^\infty \Xi_k \in {\cal F}${\rq\rq},
{\lq\lq}$X \in {\mathcal F}${\rq\rq}
and
{\lq\lq}$\Xi  \in {\cal F}\Rightarrow X \setminus \Xi \in {\cal F}${\rq\rq}.
\item[(B$_2$)]
[Countably additivity]
$F$ is a mapping from ${\cal F}$ to ${{\cal N}}$ 
satisfying:
(a):
for every $\Xi \in {\cal F}$, $F(\Xi)$ is a non-negative element in 
${{\cal N}}$
such that $0 \le F(\Xi) $
$\le I$, 
(b):
$F(\emptyset) = 0$ and 
$F(X) = I$,
where
$0$ and $I$ is the $0$-element and the identity
in ${\cal N}$
respectively.
(c):
for any countable decomposition $\{ \Xi_1,\Xi_2, \ldots \}$
of $\Xi$
$\in {\cal F}$
(i.e., $\Xi_k , \Xi \in {\cal F}$
such that
$\bigcup_{k=1}^\infty \Xi_k = \Xi$,
$\Xi_i \cap \Xi_j= \emptyset
(i \not= j)$),
it holds that
\begin{align}
&
\quad
\lim_{K \to \infty } 
{}_{{}_{{\cal N}_*}}
\langle \rho, 
F( \bigcup_{k=1}^K \Xi_k )
 \rangle_{{}_{\cal N}}
=
{}_{{}_{{\cal N}_*}}
\langle \rho, 
 F( \Xi ) 
 \rangle_{{}_{\cal N}}
\quad
(
\forall \rho \in {\frak S}^m  ({\cal N}_*)
)
\label{eq1}
\end{align}
i.e.,
$
\lim_{K \to \infty }  
F( \bigcup_{k=1}^K \Xi_k )
= F( \Xi ) 
$
in the sense of weak${}^*$ convergence in ${\cal N}$.
\end{itemize}

\par
\noindent
%
%
\par
\vskip0.3cm
\par
\par
\noindent
%

\par
\subsection{Pure Measurement Theory}
\par
\noindent

\par

Our concern in this paper is SMT in (A$_1$),
which is constructed on the base of PMT.
Thus, we begin with PMT.
\par
\rm

With any {\it system} $S$, a fundamental structure
$[{\mathcal A} \subseteq {\mathcal N} \subseteq B(H)]$
can be associated in which the 
pure
measurement theory (A$_1$) of that system can be formulated.
A {\it pure state} of the system $S$ is represented by an element
$\rho^p (\in {\frak S}^p  ({}{\cal A}^*{})$="pure state class"({\it cf.} ref.\cite{Ishi7}))
and an {\it observable} is represented by an observable 
${\mathsf{O}}{\; =} (X, {\cal F}, F)$ in ${{\cal N}}$.
Also, the {\it measurement of the observable ${\mathsf{O}}$ for the system 
$S$ with the pure state $\rho^p$}
is denoted by 
${\mathsf{M}}_{{{\cal N}}} ({\mathsf{O}}, S_{[\rho^p]})$
$\big($
or more precisely,
${\mathsf{M}}_{\cal N} ({\mathsf{O}}{\; =} (X, {\cal F}, F), S_{[\rho^p]})$
$\big)$.
An observer can obtain a measured value $x $
($\in X$) by the measurement 
${\mathsf{M}}_{\cal N} ({\mathsf{O}}, S_{[\rho^p]})$.
\par
\noindent
\par
The Axiom$^{\rm P}$ 1 presented below is 
a kind of mathematical generalization of Born's probabilistic interpretation of quantum mechanics.
\par
\noindent
{\bf{Axiom$^{\rm P}$ 1\;\;
\rm
$[$Pure Measurement$]$}}.
\it
The probability that a measured value $x$
$( \in X)$ obtained by the measurement 
${\mathsf{M}}_{{{\cal N}}} ({\mathsf{O}}$
${ \equiv} (X, {\cal F}, F),$
{}{$ S_{[\rho^p_0]})$}
%
belongs to a set 
$\Xi (\in {\cal F})$ is given by
$
\rho^p_0( F(\Xi) )
$,
if $F(\Xi)$ is essentially continuous at $\rho^p_0$
({\it cf.} ref.\cite{Ishi7}).
\rm

\par
\par
\vskip0.2cm
\par

\par
Next, we explain Axiom 2 in (A$_1$).
Let $(T,\le)$ be a tree, i.e., a partial ordered 
set such that {\lq\lq$t_1 \le t_3$ and $t_2 \le t_3$\rq\rq} implies {\lq\lq$t_1 \le t_2$ or $t_2 \le t_1$\rq\rq}\!.
In this paper,
we assume that
$T$ is finite.
Assume that
there exists an element $t_0 \in T$,
called the {\it root} of $T$,
such that
$t_0 \le t$ ($\forall t \in T$) holds.
Put $T^2_\le = \{ (t_1,t_2) \in T^2{}\;|\; t_1 \le t_2 \}$.
The family
$\{ \Phi_{t_1,t_2}{}: $
${\cal N}_{t_2} \to {\cal N}_{t_1} \}_{(t_1,t_2) \in T^2_\le}$
is called a {\it causal relation}
({\it due to the Heisenberg picture}),
\rm
if it satisfies the following conditions {}{(C$_1$) and 
(C$_2$)}.
\begin{itemize}
\item[{\rm (C$_1$)}]
With each
$t \in T$,
a fundamental structure
$[{\mathcal A}_t \subseteq {\mathcal N}_t \subseteq B(H_t)]$
is associated.
\item[{\rm (C$_2$)}]
For every $(t_1,t_2) \in T_{\le}^2$, a continuous Markov operator 
$\Phi^{t_1,t_2}{}: {\cal N}_{t_2}{\mbox{(with the weak$^*$ topology)}}
$
$ \to
$
$ {\cal N}_{t_1}$
${\mbox{(with the weak$^*$ topology)}}$ 
is defined
(i.e.,
$\Phi^{t_1,t_2} \ge 0$,
$\Phi^{t_1,t_2}(I_{{\cal N}_{t_2}})$
$
=
$
$
I_{{\cal N}_{t_1}}$
).
And it satisfies that
$\Phi^{t_1,t_2} \Phi^{t_2,t_3} = \Phi^{t_1,t_3}$ 
holds for any $(t_1,t_2)$, $(t_2,t_3) \in T_\le^2$.
\end{itemize}
\noindent
The family of pre-dual operators
$\{ \Phi^{t_1,t_2}_*{}: $
$
{\frak S}^m  (({\cal N}_{t_1})_*)
\to {\frak S}^m  (({\cal N}_{t_2})_*)
\}_{(t_1,t_2) \in T^2_\le}$
is called a
{
\it
pre-dual causal relation}
({\it
due to the Schr\"{o}dinger picture}).
%

\par
\par
\rm
Now Axiom 2 in the measurement theory (A$_1$) is presented
as follows:
\rm
\par
\noindent
{\bf{Axiom 2}
\rm[Causality]}.
\it
The causality is represented by
a causal relation 
$\{ \Phi^{t_1,t_2}{}: $
${\cal N}_{t_2} \to {\cal N}_{t_1} \}_{(t_1,t_2) \in T^2_\le}$.

\rm
\par
\par
\vskip0.2cm
\par
\noindent
\par

\par
\noindent
\par
\noindent
\par
\noindent
\vskip0.2cm
\par
\noindent
\subsection{
Linguistic Interpretation
}
\par
\noindent
\par
Next,
we have to
study the linguistic interpretation
(i.e.,
the manual of "how to use the above axioms"
)
as follows.
That is, we present the following interpretation
(D)
[=(D$_1$), (D$_2$)],
which is characterized as a kind of linguistic turn
of so-called Copenhagen interpretation
({\it cf.} refs.\textcolor{black}{\cite{Ishi7, Ishi8}}
).
That is,
\begin{itemize}
\item[(D$_1$)]
Consider the dualism composed of {\lq\lq}observer{\rq\rq} and {\lq\lq}system( =measuring object){\rq\rq}.
And therefore,
{\lq\lq}observer{\rq\rq} and {\lq\lq}system{\rq\rq}
must be absolutely separated.
\item[(D$_2$)]
Only one measurement is permitted.
And thus,
the state after a measurement
is meaningless
$\;$
since it 
can not be measured any longer.
Also, the causality should be assumed only in the side of system,
however,
a state never moves.
Thus,
the Heisenberg picture should be adopted.
And thus,
the Schr\"{o}dinger picture
is rather makeshift.
\end{itemize}
\par
\noindent
and so on.

%

\par
\noindent
\par
The following argument is a consequence of the above (D$_2$).
For each
$k=1,$
$2,\ldots,K$,
consider a measurement
${\mathsf{M}}_{{{\cal N}}} ({\mathsf{O}_k}$
${\; \equiv} (X_k, {\cal F}_k, F_k),$
$ S_{[\rho]})$.
However,
since
the (D$_2$)
says that
only one measurement is permitted,
the
measurements
$\{
{\mathsf{M}}_{{{\cal N}}} ({\mathsf{O}_k},S_{[\rho]})
\}_{k=1}^K$
should be reconsidered in what follows.
Under the commutativity condition such that
\begin{align}
&
F_i(\Xi_i) F_j(\Xi_j) 
=
F_j(\Xi_j) F_i(\Xi_i)
\label{eq2}
%
\\
&
\quad
(\forall \Xi_i \in {\cal F}_i,
\forall \Xi_j \in  {\cal F}_j , i \not= j),
\nonumber
\end{align}
we can
define the product observable
${\text{\large $\times$}}_{k=1}^K {\mathsf{O}_k}$
$=({\text{\large $\times$}}_{k=1}^K X_k ,$
$ \boxtimes_{k=1}^K {\cal F}_k,$
$ 
{\text{\large $\times$}}_{k=1}^K {F}_k)$
in ${\cal N}$ such that
\begin{align*}
({\text{\large $\times$}}_{k=1}^K {F}_k)({\text{\large $\times$}}_{k=1}^K {\Xi}_k )
=
F_1(\Xi_1) F_2(\Xi_2) \cdots F_K(\Xi_K)
\\
\;
(
\forall \Xi_k \in {\cal F}_k,
\forall k=1,\ldots,K
).
\qquad
\qquad
\nonumber
\end{align*}
Here,
$ \boxtimes_{k=1}^K {\cal F}_k$
is the smallest $\sigma$-field including
the family
$\{
{\text{\large $\times$}}_{k=1}^K \Xi_k
$
$:$
$\Xi_k \in {\cal F}_k \; k=1,2,\ldots, K \}$.
Then, 
the above
$\{
{\mathsf{M}}_{{{\cal N}}} ({\mathsf{O}_k},S_{[\rho]})
\}_{k=1}^K$
is,
under the commutativity condition (\ref{eq2}),
represented by the simultaneous measurement
${\mathsf{M}}_{{{{\cal N}}}} (
{\text{\large $\times$}}_{k=1}^K {\mathsf{O}_k}$,
$ S_{[\rho]})$.


\par
\par
Consider a tree
$(T{\; \equiv}\{t_0, t_1, \ldots, t_n \},$
$ \le )$
with the root $t_0$.
This is also characterized by
the map
$\pi: T \setminus \{t_0\} \to T$
such that
$\pi( t)= \max \{ s \in T \;|\; s < t \}$.
Let
$\{ \Phi^{t, t'} : {\cal N}_{t'}  \to {\cal N}_{t}  \}_{ (t,t')\in
T_\le^2}$
be a causal relation,
which is also represented by
$\{ \Phi_{\pi(t), t} : {\cal N}_{t}  \to {\cal N}_{\pi(t)}  \}_{ 
t \in T \setminus \{t_0\}}$.
Let an observable
${\mathsf O}_t{\; \equiv}
(X_t, {\cal F}_{t}, F_t)$ in the ${\cal N}_t$ 
be given for each $t \in T$.
Note that
$\Phi_{\pi(t), t}
{\mathsf O}_t$
$(
{\; \equiv}
(X_t, {\cal F}_{t},
\Phi_{\pi(t), t} F_t)$
)
is an observable in the ${\cal N}_{\pi(t)}$.

The pair
$[{\mathbb O}_T]
$
$=$
$[
\{{\mathsf O}_t \}_{t \in T}$,
$\{ \Phi^{t, t'} : {\cal N}_{t'}  \to {\cal N}_{t}  \}_{ (t,t')\in
T_\le^2}$
$]$
is called a
{\it sequential causal observable}.
%
%
For each $s \in T$,
put $T_s =\{ t \in T \;|\; t \ge s\}$.
And define the observable
${\widehat{\mathsf O}}_s
\equiv ({\text{\large $\times$}}_{t \in T_s}X_t, \boxtimes_{t \in T_s}{\cal F}_t, {\widehat{F}}_s)$
in ${\cal N}_s$
as follows:
\par
\noindent
\begin{align}
\widehat{\mathsf O}_s
&=
\left\{\begin{array}{ll}
{\mathsf O}_s
\quad
&
\!\!\!\!\!\!\!\!\!\!\!\!\!\!\!\!\!\!
\text{(if $s \in T \setminus \pi (T) \;${})}
\\
{\mathsf O}_s
{\text{\large $\times$}}
({}\bigtimes_{t \in \pi^{-1} ({}\{ s \}{})} \Phi_{ \pi(t), t} \widehat {\mathsf O}_t{})
\quad
&
\!\!\!\!\!\!
\text{(if $ s \in \pi (T) ${})}
\end{array}\right.
\label{eq3}
\end{align}
if
the commutativity condition holds
(i.e.,
if the product observable
${\mathsf O}_s
{\text{\large $\times$}}
({}\bigtimes_{t \in \pi^{-1} ({}\{ s \}{})} \Phi_{ \pi(t), t}
$
$\widehat {\mathsf O}_t{})$
exists)
for each $s \in \pi(T)$.
Using (\ref{eq3}) iteratively,
we can finally obtain the observable
$\widehat{\mathsf O}_{t_0}$
in ${\cal N}_{t_0}$.
The
$\widehat{\mathsf O}_{t_0}$
is called the realization
(or,
realized causal observable)
of
$[{\mathbb O}_T]$.
%

\noindent

%
%
\noindent
\noindent
\subsection{Statistical Measurement Theory
}

\par
\noindent
\par
Let $[{\mathcal A} \subseteq {\mathcal N} \subseteq B(H)]$
be a fundamental structure.
We shall introduce the following notation:
\rm
It is usual to consider that
we do not know the pure state
$\rho_0^p$
$(
\in
{\frak S}^p  ({}{\cal A}^*{})
)$
when
we take a measurement
${\mathsf{M}}_{{{\cal N}}} ({\mathsf{O}}, S_{[\rho_0^p]})$.
That is because
we usually take a measurement ${\mathsf{M}}_{{{\cal N}}} ({\mathsf{O}},
S_{[\rho_0^p]})$
in order to know the state $\rho_0^p$.
Thus,
when we want to emphasize that
we do not know the state $\rho_0^p$,
${\mathsf{M}}_{{{\cal N}}} ({\mathsf{O}}, S_{[\rho_0^p]})$
is denoted by
${\mathsf{M}}_{{{\cal N}}} ({\mathsf{O}}, S_{[\ast]})$.
Also,
when we know the distribution $\rho_0^m$
$( \in {\frak S}^m({\cal N}_*) )$
of the unknown state
$\rho_0^p$,
the
${\mathsf{M}}_{{{\cal N}}} ({\mathsf{O}}, S_{[\rho_0^p]})$
is denoted by
${\mathsf{M}}_{{{\cal N}}} ({\mathsf{O}}, S_{[\ast]}
({ \rho_0^m }) )$.
The $\rho_0^m$
is called a mixed state.
In Bayesian statistics, the mixed state $\rho_0^m$
may be called a "subjective pretest state"
({\it cf.} refs.\textcolor{black}{\cite{Ishi7}--\cite{Ishi10}}
).

\par
\vskip0.3cm

\par

\par
The Axiom$^{\rm S}$ 1 presented below 
is also
a kind of mathematical generalization of Born's probabilistic interpretation of quantum mechanics.

\par

\par
\noindent
{\bf{Axiom$^{\rm S}$\;1\;
\rm
\;[Statistical measurement]}}.
\it
The probability that a measured value $x$
$( \in X)$ obtained by the measurement 
${\mathsf{M}}_{{{\cal N}}} ({\mathsf{O}}$
${ \equiv} (X, {\cal F}, F),$
{}{$ S_{[\ast]}({ \rho_0^m }) )$}
%
belongs to a set 
$\Xi (\in {\cal F})$ is given by
$
\rho_0^m ( F(\Xi) )
$
$($
$=
{}_{{{\cal N}_*}}\langle
\rho_0^m,
F(\Xi)
\rangle_{{\cal N}}$
$)$.
\rm

\subsection{Examples of SMT
(Bertrand paradox)}
\par
\rm
\par
\par
\noindent
\par
Consider classical systems in a commutative $W^*$-algebra
$L^\infty ( \Omega, m )$.
We can define the exact observable
${\mathsf O}_E
=( \Omega, {\mathcal B}_\Omega, F_E )$ in $L^\infty ( \Omega, m )$
such that
$$
[F_E(\Xi)](\omega )=
\chi_{{}_\Xi}
(\omega)
=
\begin{cases}
1 \quad & (\omega \in \Xi )
\\
0 \quad & (\omega \notin \Xi )
\end{cases}
\\
\qquad
(\forall \omega \in \Omega, \;\;
\Xi \in {\mathcal B}_\Omega )
$$
Here, we have the following problem:
\begin{itemize}
\item[(E$_1$)]
Can the measurement
${\mathsf M}_{L^\infty ( \Omega, m )}({\mathsf O}_E, S_{[\ast]}(\rho) )$
that represents "at random"
be determined uniquely?
\end{itemize}
This question is of course denied by so-called Bertrand paradox.
Here,
let us review the argument about the Bertrand paradox
({}{\it cf.} ${{{}}}$\cite{Isaa, Ishi6}{}).
Consider the following problem:
\begin{itemize}
\item[(E$_2$)]
Given a circle with the radius 1. 
Suppose a chord of the circle is chosen at random. What is the probability that the chord is shorter than 
${\sqrt 3 }$?
%
%
%
\end{itemize}
\par
\noindent
This problem is not well-posed.
That is because this (E$_2$) is essentially the same as the following
problem (E$_3$): 
\par
\noindent
\par
\noindent
\unitlength=0.8mm
\begin{picture}(200,75)
\put(40,0){
\put(53,67){(Fig.0)}
\put(64,53){\scriptsize $\l$}
\bezier{200}(42,55)(56,56))(63,54)
\bezier{200}(70,53)(78,51))(85,45)
\allinethickness{0.5mm}
\put(42,55){\line(4,-1){44}}
\allinethickness{0.7mm}
\put(60,30){\circle{60}}
}
\end{picture}
\par
\noindent
\par
\noindent
\begin{itemize}
\item[(E$_3$)]
In Fig. 0,
does there exist the most natural coordinate by which 
a chord $\l$ is represented?
Or, is the most natural coordinate determined uniquely?
\end{itemize}
This is of course denied as follows.
\par
\par
\noindent
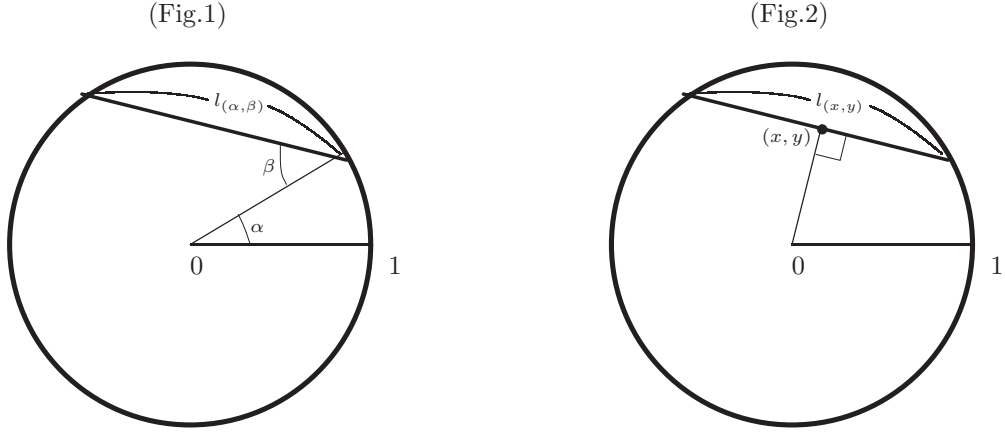
\unitlength=0.8mm
\begin{picture}(200,75)
\put(70,32){\scriptsize $\alpha$}
\spline(70,30)(69,33)(68,35)
\put(72,42){\scriptsize $\beta$}
\spline(76,40)(75,42)(75,47)
\put(153,67){({Fig.$2$})}
\put(53,67){({Fig.$1$})}
\put(155,47){\scriptsize $(x,y)$}
\put(164,48){$\bullet$}
\put(160,25){$0$}
\put(193,25){$1$}
\put(60,25){$0$}
\put(93,25){$1$}
\put(0,0){
\put(64,53){\scriptsize $\l_{(\alpha, \beta{})}$}
\bezier{200}(42,55)(56,56))(63,54)
\bezier{200}(73,53)(78,51))(85,45)
}
\put(100,0){
\put(64,53){\scriptsize $\l_{(x,y{})}$}
\bezier{200}(42,55)(56,56))(63,54)
\bezier{200}(73,53)(78,51))(85,45)
\path(64,45)(68,44)(69,48)}
\allinethickness{0.2mm}
\put(60,30){\line(1,0){30}}
\put(160,30){\line(1,0){30}}
\put(60,30){\line(5,3){26}}
\put(160,30){\line(1,4){5}}
\allinethickness{0.05mm}
\allinethickness{0.5mm}
\put(42,55){\line(4,-1){44}}
\put(142,55){\line(4,-1){44}}
\allinethickness{0.7mm}
\put(60,30){\circle{60}}
\put(160,30){\circle{60}}
\end{picture}
\vskip1.0cm

\par
\noindent
\bf
[The first answer ({}{Fig.$1$})].
\rm
In {Fig.$1$}, we see that
the chord $\l$ is represented by a point
$({}\alpha, \beta{}) $
in the rectangle
$\Omega_1$
$\equiv$
$\{ ({}\alpha, \beta{}) \; | \; 0 < \alpha \le 2 \pi, \;
0 < \beta \le \pi/2 \text{({}radian{})}\}$.
That is, we have the following identification:
\begin{align*}
\Omega
(
=
\mbox{the set of all chords}) \ni \l_{({}\alpha, \beta{})}
\approx
({}\alpha , \beta{}) \in \Omega_1
( \subset {\mathbb R}^2 ).
\end{align*}
Under the identification,
we get the natural probability measure
$m_1$
on
$\Omega_1$
such that
$m_1 ({}A{}) = \frac{ \text{Meas}[{}A]}{\text{Meas}[{}\Omega_1{}] }
= \frac{ \text{Meas} [{}A]}{ \pi^2 }$
$({}\forall A \in {\cal B}_{\Omega_1}{})$,
where
{\lq\lq} Meas{\rq\rq} = {\lq\lq} Lebesgue measure{\rq\rq}$\!\!\!\!.\; \;$
Therefore,
we have a natural measurement
${\mathsf M}_{L^\infty ( \Omega_1, m_1 )}({\mathsf O}_E \equiv
(\Omega_1, {\mathcal B}_{\Omega_1} , F_E), S_{[\ast]}(1) )$.
Put
$$
\Xi^1_{\sqrt{3}}=\{ (\alpha , \beta )\in
\Omega_1
\;:\;
\text{"the length of } \l_{({}\alpha , \beta{})}" < {\sqrt 3}
\}
$$
Then, Axiom$^{\mbox{\scriptsize S}}$ 1 says that
the probability that a measured value $(\alpha, \beta )$
belongs to
$\Xi^1_{\sqrt{3}}$
is given by
\par
\noindent
\begin{align*}
&
\int_{\Omega_1} [F_E (\Xi^1_{\sqrt{3}})](\omega_1) \; m_1( d \omega_1 )=
\int_{\Xi^1_{\sqrt{3}}}\;\; 1 \;\; m_1( d \omega_1 )
\\
=
&
m_1
(\{ \l_{({}\alpha , \beta{})} \approx (\alpha, \beta) \in \Omega_1 \; | \;
\text{ "the length of } \l_{({}\alpha , \beta{})}" \le {\sqrt 3} \}{})
\\
=
&
\frac{\text{Meas}[{}\{ ({}\alpha, \beta{}) \; | \; 0 \le \alpha \le 2 \pi,
\;\pi/6 \le \beta \le \pi/2 \}]}
{\text{Meas}[{}\{ ({}\alpha, \beta{}) \; | \; 0 \le \alpha \le 2 \pi,
\;0 \le \beta \le \pi/2\}]}
\\
=
&
\frac{2 \pi \times ( \pi/3)}
{\pi^2}
= \frac{2}{3}.
\end{align*}
\par
\noindent
\bf
[The second answer ({}{Fig.$2$})].
\rm
In {Fig.$2$}, we see that
the chord $\l$ is represented by a point
$(x,y) $
in the circle
$\Omega_2$
$\equiv$
$\{ (x,y) \; | \; 
x^2 + y^2 < 1
\}$.
That is, we have the following identification:
\begin{align*}
\Omega
(
=
\mbox{the set of all chords})
\ni \l_{(x,y)}
\approx
(x,y) \in \Omega_2
(
\subset {\mathbb R}^2 ).
\end{align*}
Under the identification,
we get the natural probability measure
$m_2$
on $\Omega_2$
such that
$m_2 ({}A{}) = \frac{ \text{Meas}[{}A]}{\text{Meas}[{}\Omega_2{}] }
= \frac{ \text{Meas} [{}A]}{  \pi }$
$({}\forall A \in {\cal B}_{\Omega_2}{})$.
Therefore,
we have a natural measurement
${\mathsf M}_{L^\infty ( \Omega_2, m_2 )}({\mathsf O}_E \equiv
(\Omega_2, {\mathcal B}_{\Omega_2} , F_E), S_{[\ast]}(1) )$.
Put
$$
\Xi^2_{\sqrt{3}}=\{ (x,y )\in
\Omega_2
\;:\;
\text{"the length of } \l_{({}x,y{})}" < {\sqrt 3}
\}
$$
Then, Axiom$^{\mbox{\scriptsize S}}$ 1 says that
the probability that a measured value $(x,y )$
belongs to
$\Xi^2_{\sqrt{3}}$
is given by
\par
\noindent
\begin{align*}
&
\int_{\Omega_2} [F_E (\Xi^2_{\sqrt{3}})](\omega_2)  \; m_2( d \omega_2 )=
\int_{\Xi^2_{\sqrt{3}}}\;\; 1 \;\; m_2( d \omega_2 )
\\
=
&
m_2
(\{ \l_{({}x,y)} \approx (x,y)  \in \Omega_2 \; | \;
\text{ "the length of } \l_{(x,y)}" \le {\sqrt 3} \}{})
\\
=
&
\frac{\text{Meas} [{}
\{ (x,y) \; | \; 1/4 \le x^2+y^2 \le 1 \}
]}
{\pi}
= \frac{3}{4}.
\end{align*}
%
%
%
\par
Therefore, the first answer and the second answer say that
the question (E$_1$) is denied. Also note that this is not a paradox since quantum language is a language.
And thus,
it is certain that, for each $i=1,2$, there exist some phenomena such that they are represented by the measurement
${\mathsf M}_{L^\infty ( \Omega_i, m_i )}({\mathsf O}_E, S_{[\ast]}(1) )$.
However, it should be noted that,
if the state space $\Omega$ is finite,
the concept of "at random"
(or, "having no information")
is meaningful in the sense of ref. \cite{Ishi10}.


\section{
Bayes Method in classical
$L^\infty(\Omega, m)$
}
\subsection{
Notations
}

\par
\par
\noindent
\par
As mentioned in the previous section, in this paper, we devote ourselves to
$$
\fbox{SMT${}^{W^*}$}=\overset{\roman measurement}{\fbox{Axiom$^{\roman S}$ 1}}+\overset{\roman causality}{\fbox{Axiom 2}}+\overset{\text{ interpretation(= how to use axioms)}}{\fbox{Linguistic Interpretation}}
$$
in classical sysytems.
Thus, the fundamental structure $[{\mathcal A} \subseteq {\mathcal N} \subseteq B(H)]$
is represented by
$$
[
C_0(\Omega ) \subseteq L^\infty(\Omega, m ) \subseteq B(L^2(\Omega, m ))
]
$$
where
the measure space $(\Omega, {\mathcal B}_\Omega , m )$ is assumed to satisfy that
$$
m(\{\omega \})=0 \;\;(\forall \omega \in \Omega),
\qquad
m(D ) > 0 \;\; (\forall \mbox{open set $D ( \subseteq \Omega \mbox{: locally compact space})$)}
$$
also,
the Borel $\sigma$-field
${\mathcal B}_\Omega$
is defined as
the smallest $\sigma$-field including
all open sets in $\Omega$.
It should be noted that
the dual Banach space
of
$L^1(\Omega )$
(i.e.,
the space composed of all integral functions on $\Omega$)
is equal to $L^\infty (\Omega )$,
that is,
$L^\infty (\Omega )_*
=L^1(\Omega )$.
Therefore, the bi-linear functional
$\rho ( f)$ is represented by
\begin{align*}
&
\rho ( f)=
{}_{{}_{L^1(\Omega)}} \langle  \rho, F \rangle {}_{{}_{L^\infty(\Omega)}}
=
\int_\Omega f(\omega ) \rho (\omega ) m( d \omega )
\\
(\forall
f \in L^\infty (\Omega, m ) &
\mbox{(or, in short $L^\infty(\Omega)$)}, \;\;
\rho \in L^1(\Omega,m )\mbox{(or, in short $L^1(\Omega)$)})
\end{align*}
({\it cf.$\;$}\textcolor{black}{\cite{Yosi}}).
Also, note that
the mixed state space
${\frak S}^m({\mathcal N}_* )$ is characterized such as
$$
{\frak S}^m({\mathcal N}_* )=
L^1_{+1}(\Omega)
\equiv
\{
\rho \in L^1(\Omega )
\;:\;
\rho \ge0,
\;\;
|| \rho||_{L^1(\Omega )}
\equiv
\int_\Omega |\rho(\omega )| m (d \omega )
=
1
\}.
$$

\par

\subsection{
Bayes Method
in Classical
$L^\infty(\Omega, m)$
}

\par
\noindent

\rm 
Let ${\mathsf O} \equiv (X, {\cal F}, F)$ be an observable in a commutative $W^*$-algebra $L^\infty(\Omega, m)$.
And let
${\mathsf O}' \equiv (Y, {\cal G}, G)$ be any observable in $L^\infty(\Omega, m)$.
Consider the product observable
${\mathsf O} \times {\mathsf O}' \equiv (X\times Y, {\cal F}
\bigstimes
{\cal G}, F \times G)$ in $L^\infty(\Omega, m)$.

Assume that
we know that
the measured value $(x,y)$ obtained by a simultaneous measurement
${\mathsf M}_{L^\infty(\Omega, m)}(
{\mathsf O} \times {\mathsf O}',
S_{[*]}
{({\rho_0})}
)$
belongs to
$\Xi \times Y \;(\in {\cal F} \boxtimes {\cal G} )$.
Then,
by Axiom${}^{{\rm S}}$ 1,
we say that
\begin{itemize}
\item[{\rm (F$_1$)}]
the probability $
P_\Xi (G(\Gamma))
$
that
$y$ belongs to $\Gamma (\in {\cal G})$
is given by
$$
\!\!\!
P_\Xi (G(\Gamma))
= \frac{\int_{\Omega} [F(\Xi) \cdot G(\Gamma)](\omega)
\;
\rho_0(\omega)
\;
m (d \omega) }{
\int_{\Omega} [F(\Xi)](\omega)
\;
\rho_0(\omega)
\;
m(d \omega )}
\;\;
(\forall \Gamma \in {\cal G}).
$$
\end{itemize}
Thus,
we can assert that:
\begin{itemize}
\item[(F$_2$)]
When
we
know that
a measured value
obtained by
a
measurement
${\mathsf M}_{L^\infty(\Omega, m)}(
{\mathsf O} \equiv (X, {\cal F}, F)
, S_{[*]}
{(\rho_0 )}
)$
belongs to
$\Xi$,
there is a reason to
infer
that
the mixed state after the measurement
is
equal to
$\rho_0^a$
$( \in L^1_{+1} (\Omega ))$,
where
$$
\rho_0^a( \omega)= \frac{ [F(\Xi_)](\omega)
\;
\rho_0 (\omega) }{
\int_{\Omega} [F(\Xi)](\omega)
\;
\rho_0 (\omega ) 
\;
m(d \omega )}
\quad
(\forall \omega \in \Omega ).
$$
\end{itemize}
After all, 
we can define the Bayes operator
$[B_{{\mathsf{O}} }^0(\Xi)]:
L^1_{+1}(\Omega)$
$
\to L^1_{+1}( \Omega)$
such that
\begin{itemize}
\item[(F$_3$)]
$
\quad
\overset{\text{(pretest state)}}{
\underset{
(\in L^1_{+1}(\Omega))}
{
\fbox{$
{\rho_0}
$
}
}
}
\xrightarrow[\text{\tiny Bayes operator}]{\mbox{{
$\qquad [B_{{\mathsf{O}} }^0(\Xi)]
 \qquad$}}}
\overset{\text{(posttest state)}}{
\underset{(
\in L^1_{+1}(\Omega))}
{\fbox{$
{\rho_0^a}
$
}}
}
$
\end{itemize}
%

\rm
\par
\vskip0.3cm
\par
\noindent
{\it Remark 1.}
The above (F$_2$) superficially contradicts the linguistic interpretation (D$_2$),
which says that
"a state never move".
In this sense,
the above (F$_2$) is convenient and makeshift.
For the precise argument, see \cite{Ishi7, Ishi8}.
That is, in spite of the linguistic interpretation (D$_2$),
we admit
the wavefunction collapse such as (F$_2$).

\rm


\par

\par
\par
\noindent

\subsection{
Bayes-Kalman Method
in Classical
$L^\infty(\Omega, m)$
}

\par
\noindent

Let
$t_0$
be the root
of a tree $T$.
Let
$[{\mathbb O}_T{}]$
$=$
$[{}
\{ {\mathsf{O}}_t  ({}\equiv ({}X_t   ,$
$ {\cal F}_{t} , {F}_t  ))
\}_{ t \in T} ,
\{  \Phi^{t_1,t_2}{}: $
$L^\infty (\Omega_{t_2}) \to L^\infty (\Omega_{t_1}) \}_{(t_1,t_2) \in T^2_\le }$
$]$
be
a sequential causal observable
with
the realization
$\widehat{\mathsf{O}}_{t_0} $
$\equiv$
$(\bigtimes_{t \in T } X_t  , $
$\bigstimes_{t \in T } {\cal G}_t,$
${\widehat F}_{t_0})$
in $L^\infty (\Omega_{t_0})$.
And let
$[{\mathbb O}_T^\times{}]$
$=$
$[{}
\{ {\mathsf{O}}_t^\times  ({}\equiv ({}X_t \times Y_t  ,$
$ {\cal F}_{t} \boxtimes 
{\cal G}_t , {F}_t \times G_t ))
\}_{ t \in T} ,
\{  \Phi^{t_1,t_2}{}: $
$L^\infty (\Omega_{t_2}) \to L^\infty (\Omega_{t_1}) \}_{(t_1,t_2) \in T^2_\le }$
$]$
be
a sequential causal observable
with
the realization
$\widehat{\mathsf{O}}_{t_0}^\times $
$\equiv$
$(\bigtimes_{t \in T } (X_t \times Y_t) , $
$\bigstimes_{t \in T }
({\cal F}_{t} \boxtimes {\cal G}_t),$
${\widehat H}_{t_0})$
in $L^\infty (\Omega_{t_0})$.
Thus we have the statistical measurement
${\mathsf M}_{C(\Omega_{t_0})}
(\widehat{\mathsf{O}}_{t_0}^\times, S_{[\ast]}(\rho_0 ) )$,
where
$\rho_0 \in L^1_{+1} (\Omega_{t_0} )$.
Assume that
we know that
the measured value $(x,y)$
$(
=
(
(x_t)_{t\in T},
(y_t)_{t\in T},
)
\in (\bigtimes_{t \in T} X_t)
\bigtimes
(\bigtimes_{t \in T} Y_t)
)$
obtained by the measurement
${\mathsf M}_{C(\Omega_{t_0})}
(\widehat{\mathsf{O}}_{t_0}^\times, S_{[\ast]}(\rho_0 ) )$
belongs to
$(\bigtimes_{t \in T} \Xi_t)
\bigtimes 
$
$
(\bigtimes_{t \in T}Y_t) \;$
$(\in 
(\boxtimes_{t \in T}{\mathcal F}_t)
\boxtimes 
(\boxtimes_{t \in T} {\mathcal G}_t)
)$.
Then,
by Axiom${}^{{\rm S}}$ 1,
we can infer that
\begin{itemize}
\item[{\rm (G$_1$)}]
the probability $
P_{\times_{t \in T} \Xi_t} (
(G_t (\Gamma_t))_{t \in T}
)
$
that
$y$ belongs to $\bigtimes_{t \in T} \Gamma_t (\in \boxtimes_{t \in T} {\cal G}_t)$
is given by
\begin{align}
&
P_{\times_{t \in T} \Xi_t} (
(G_t (\Gamma_t))_{t \in T}
)
\nonumber
\\
= 
&
\frac{\int_{\Omega_0} 
[
{\widehat H}_{t_0}
(
(\bigtimes_{t \in T} \Xi_t)
\bigtimes 
(\bigtimes_{t \in T} \Gamma_t) 
)
]
(\omega_0)
\;
\rho_0 (\omega_0 ) 
\;
m_0 (d \omega_0) }{
\int_{\Omega_0} [
{\widehat H}_{t_0}
(\bigtimes_{t \in T} \Xi_t)
\bigtimes 
(\bigtimes_{t \in T}Y_t) 
](\omega_0)
\;
\rho_0 (\omega_0 ) 
\;
m_0(d \omega_0 )}
\label{eq4}
\\
&
\quad
(\forall \Gamma_t \in {\cal G}_t, t \in T).
\nonumber
\end{align}
\end{itemize}
Let $s \in T$ be fixed.
Assume that
$$
\Gamma_t = Y_t
\quad (\forall t \in T \mbox{ such that $t \not= s$})
$$
Thus, putting
${\widehat P}_{\times_{t \in T} \Xi_t} (
G_s (\Gamma_s)
)=P_{\times_{t \in T} \Xi_t} (
(G_t (\Gamma_t))_{t \in T}
)$,
we see that
$
{\widehat P}_{\times_{t \in T} \Xi_t} \in L^1_{+1}(\Omega_s, m_s )
$.
%
%
%
%
That is,
there uniquely exists 
$
\rho_s^a
\in L^1_{+1} (\Omega_s , m_s)
$
such that
\begin{align}
&
{\widehat P}_{\times_{t \in T} \Xi_t} (
(G_s (\Gamma_s))
=
{}_{{}_{L^1(\Omega_s)}} \langle  \rho_s^a, G_s (\Gamma_s) \rangle {}_{{}_{L^\infty(\Omega_s)}}
=
\int_{\Omega_s} [G_s (\Gamma_s)](\omega_s ) \rho_s^a(\omega_s )
m_s ( d \omega_s )
\label{eq5}
\end{align}
for any observable
$(Y_s, {\cal G}_s, G_s)$
in
$L^\infty (\Omega_s)$.
That is because
the linear functional
${\widehat P}_{\times_{t \in T} \Xi_t} 
: L^\infty(\Omega_s ) \to 
{\mathbb C}
$
(complex numbers)
is weak$^*$ continuous.
After all, 
we can define the Bayes-Kalman operator
$[B_{\widehat{\mathsf{O}}_{t_0} }^s(\bigtimes_{t \in T} \Xi_t)]:
L^1_{+1}(\Omega_{t_0})$
$
\to L^1_{+1}( \Omega_s)$
such that
\begin{itemize}
\item[(G$_2$)]
$
\quad
\overset{\text{(pretest state)}}{
\underset{
(\in L^1_{+1}(\Omega_{t_0}))}
{\fbox{$
{\rho_0}
$
}}
}
\xrightarrow[\text{\tiny Bayes-Kalman operator}]{\mbox{
$\qquad [B_{\widehat{\mathsf{O}}_{t_0} }^s(\bigtimes_{t \in T} \Xi_t)] \qquad$}}
\overset{\text{(posttest state)}}{
\underset{(
\in L^1_{+1}(\Omega_s))}
{\fbox{$
{\rho_0^a}
$
}}
}
$
\end{itemize}
which is the generalization of the \textcolor{black}{(F$_3$)}
({\it cf}. \cite{Kal1, Ishi5, Kiku1}).
Here, the ref. \cite{Kal1} is of course great. However,
now we think that our arguments in refs. \cite{Kiku1, Ishi6} were too abstract and philosophical, that is, we should have presented concrete examples much more.
This will be done in the following sections.
%


\par
\rm
\vskip1.0cm
\rm

\section{A simle example of Kalman filter}
\subsection{The computable situation}
\par
\noindent
\par
Put $T=\{0,1,2, \cdots, n\}$,
and define
the parents map
$\pi : T\setminus \{0\} \to T$
such that
$\pi(k) = k-1$.
For each $k \in T$, define the commutative $W^*$-algebra: 
\begin{align*}
{\mathcal N}_k
=L^\infty ( \Omega_k , m_k )=
L^\infty ( {\mathbb R}, d \omega )
\quad
\mbox{(where $d \omega$ is the usual Lebesgue measure on ${\mathbb R}$)}
\end{align*}
and thus, the pre-dual Banach space $({\mathcal N}_k)_*$ is defined as follows.
\begin{align*}
({\mathcal N}_k)_*
=L^1 ( \Omega_k , m_k )=
L^1 ( {\mathbb R}. d \omega )
\quad
(k \in T =\{0,1,2, \cdots, n \})
\end{align*}
Consider the sequential observable
$[{\mathbb O}_T]
$
$=$
$[
\{{\mathsf O}_t \}_{t \in T}$,
$\{ \Phi^{t-1, t} : {\cal N}_{t}  \to {\cal N}_{t-1}  \}_{T=1,2, \cdots, n}
$
$]$,
and assume the initial state
$\rho_0 \in L^1_{+1} (\Omega_0 , m_0 )$.
\par
\noindent
Thus, we have the following situation:
{
\small
\begin{align*}
\overset{\mbox{initial state $\rho_0$}}{\underset{{\mathsf O}_0=(X_0, {\mathcal F}_0 F_0)}{\fbox{\mbox{$L^\infty(\Omega_0, m_0 )$}}}}
\xleftarrow[]{\Phi^{0,1}}
\underset{{\mathsf O}_1=(X_1, {\mathcal F}_1 F_1)}{\fbox{\mbox{$L^\infty(\Omega_1, m_1 )$}}}
\xleftarrow[]{\Phi^{1,2}}
\cdots
\xleftarrow[]{\Phi^{s-1,s}}
\underset{{\mathsf O}_s=(X_s, {\mathcal F}_s F_s)}{\fbox{\mbox{$L^\infty(\Omega_s, m_s )$}}}
\xleftarrow[]{\Phi^{s,s+1}}
\cdots
\xleftarrow[]{\Phi^{n-1,n}}
\underset{{\mathsf O}_n=(X_n, {\mathcal F}_n F_n)}{\fbox{\mbox{$L^\infty(\Omega_n, m_n )$}}}
\end{align*}
}
or, equivalently,
{
\small
\begin{align*}
\overset{\mbox{initial state $\rho_0$}}{\underset{{\mathsf O}_0=(X_0, {\mathcal F}_0, F_0)}{\fbox{\mbox{$L^1(\Omega_0, m_0 )$}}}}
\xrightarrow[]{\Phi_*^{0,1}}
\underset{{\mathsf O}_1=(X_1, {\mathcal F}_1, F_1)}{\fbox{\mbox{$L^1(\Omega_1, m_1 )$}}}
\xrightarrow[]{\Phi_*^{1,2}}
\cdots
\xrightarrow[]{\Phi_*^{s-1,s}}
\underset{{\mathsf O}_s=(X_s, {\mathcal F}_s, F_s)}{\fbox{\mbox{$L^1(\Omega_s, m_s )$}}}
\xrightarrow[]{\Phi_*^{s,s+1}}
\cdots
\xrightarrow[]{\Phi_*^{n-1,n}}
\underset{{\mathsf O}_n=(X_n, {\mathcal F}_n, F_n)}{\fbox{\mbox{$L^1(\Omega_n, m_n )$}}}
\end{align*}
}
Here the initial state
$\rho_0 (\in L^1_{+1} (\Omega_0 , m_0 ))$
is defined by
\begin{align}
\rho_0(\omega_0)
=
\frac{1}{\sqrt{2  \pi} \sigma_0}
\exp[-
\frac{(\omega_0-\mu_0)^2}{2 \sigma_0^2}
]
\qquad( \forall \omega_0 \in \Omega_0)
\label{eq6}
\end{align}
where it is assumed that 
$\mu_0$ and $\sigma_0$ are known.
\par
Also, for each $t \in T =\{0,1, \cdots, n \}$, consider the observable
${\mathsf O}_t=(X_t, {\mathcal F}_t, F_t)$ 
$=({\mathbb R}, {\mathcal B}_{{\mathbb R}}, F_t )$
in $L^\infty ( \Omega_t, m_t )$
such that
\begin{align}
[F_t (\Xi_t )](\omega_t)
=
\int_{\Xi_t}
\frac{1}{\sqrt{2 \pi} {q}_t}
\exp[
-\frac{(x_t -c_t \omega_t -d_t)^2}{2{q}_t^2}
]
dx_t
\equiv
\int_{\Xi_t}
f_{x_t}(\omega_t) dx_t
\quad(\forall \Xi_t \in {\mathcal F}_t, \;\; \forall \omega_t \in \Omega_t )
\label{eq7}
\end{align}
where it is assumed that 
$c_t$, $d_t$ and $q_t$ are known
$(t \in T)$.
\par
And further, the causal operator
$\Phi^{t-1.t}:L^\infty (\Omega_t) \to L^\infty ( \Omega_{t-1})$
is defined by
\begin{align}
&
[\Phi^{t-1,t} \widetilde{f}_{x_{t}}](\omega_{t-1} )
=
\int_{- \infty}^{\infty}
\frac{1}{\sqrt{2 \pi} r_t}
\exp[
-\frac{(\omega_t - a_t \omega_{t-1} -b_t )^2}{2r_t^2}
]
\widetilde{f}_{x_{t}}
)
d \omega_{t}
\equiv
f_{t-1}(\omega_{t-1})
\label{eq8}
\\
&
\qquad \qquad \qquad
(\forall \widetilde{f}_{x_{t}} \in L^\infty(\Omega_{t}, m_{t}),
\;\;
\forall \omega_{t-1} \in \Omega_{t-1})
\nonumber
\end{align}
where it is assumed that 
$a_t$, $b_t$ and $r_t$ are known
$(t \in T)$.
\par

%
%
%
%
%
%
%

\par
\noindent
Or, equivalently,
the pre-dual causal operator
$\Phi_*^{t-1.t}:L^1_{+1} (\Omega_{t-1}) \to L^1_{+1} ( \Omega_{t})$
is defined by
\begin{align}
&
[\Phi^{t-1,t}_* \widetilde{\rho}_{t-1}](\omega_t )
=
\int_{- \infty}^{\infty}
\frac{1}{\sqrt{2 \pi} r_t}
\exp[
-\frac{(\omega_t - a_t \omega_{t-1} -b_t )^2}{2r_t^2}
]
\widetilde{\rho}_{t-1}(\omega_{t-1})
d\omega_{t-1}
\label{eq9}
\\
&
\qquad \qquad \qquad
(\forall \widetilde{\rho}_{t-1} \in L^1_{+1}(\Omega_{t-1}, m_{t-1}),
\;
\forall \omega_t \in \Omega_t)
\nonumber
\end{align}

Now we have
the sequential observable
$[{\mathbb O}_T]
$
$=$
$[
\{{\mathsf O}_t \}_{t \in T}$,
$\{ \Phi^{t-1, t} : {\cal N}_{t}  \to {\cal N}_{t-1}  \}_{T=1,2, \cdots, n}
$.
Let
$\widehat{\mathsf O}_{0}$
$(\bigtimes_{t=0}^n X_t,
\boxtimes_{t=0}^n {\mathcal F}_t,
{\widehat F}
)$
be its realization.
Then we have the following problem:
\begin{itemize}
\item[(H)]
Assume that
a measured value
$(x_0, x_2, \cdots, x_n )$
$(\in \bigtimes_{t=0}^n X_t)$
is obtained by the measurement
${\mathsf M}_{L^\infty (\Omega_0)}$
$(\widehat{\mathsf O}_{0},$
$
S_{[\ast]}(\rho_0 )
)$.
Let $s(\in T)$ be fixed. Then,
calculate
the Bayes-Kalman operator
$[B_{\widehat{\mathsf{O}}_{0} }^s(\bigtimes_{t \in T} \{x_t\})]
( \rho_0 )$
in (G$_2$), where
$$
[B_{\widehat{\mathsf{O}}_{0} }^s(\bigtimes_{t \in T} \{x_t\})]
( \rho_0 )
=
\lim_{\Xi_t \to x_t \;(t\in T)}
[B_{\widehat{\mathsf{O}}_{0} }^s(\bigtimes_{t \in T} \Xi_t)](\rho_0)
$$
\end{itemize}

\subsection{
Bayes-Kalman operator
$[B_{\widehat{\mathsf{O}}_{0} }^s (\bigtimes_{t \in T} \{x_t \})]:
L^1_{+1}(\Omega_{s})$
$
\to L^1_{+1}( \Omega_n)$
}

In what follows, we solve the problem (H).
For this,
it suffices to find the
$\rho_s \in L^1_{+1} (\Omega_s )$ such that
\begin{align}
\lim_{ \Xi_t \to x_t \;\;(t\in T)}
\frac{\int_{\Omega_0} [{\widehat F}_{0}((\bigtimes_{t=0}^n \Xi_{t}) \times \Gamma_s)](\omega_0)  \;\;
\rho_0 (\omega_0 ) d\omega_0}{\int_{\Omega_0}  [{\widehat F}_{0}(\bigtimes_{t=0}^n \Xi_{t}) ](\omega_0)\;\;
\rho_0 (\omega_0 )   d\omega_0}
=
\int_{\Omega_s} [G_s(\Gamma_s )](\omega_s) 
\;\;
\rho_s (\omega_s ) d \omega_s
\quad (\forall \Gamma_s \in {\mathcal F}_s )
\nonumber
\end{align}
Let us calculate
$\rho_s= [B_{\widehat{\mathsf{O}}_{0} }^s (\bigtimes_{t \in T} \{x_t \})]({\rho}_0)$
as follows.

\begin{align}
&
\int_{\Omega_0} [{\widehat F}_{0}((\bigtimes_{t=0}^n \Xi_{t}) \times \Gamma_s)](\omega_0)  \;\;
\rho_0 (\omega_0 ) d\omega_0
\nonumber \\
=
&
{}_{{}_{L^1(\Omega_0)}} \langle  \rho_0, 
{\widehat F}_{0}((\bigtimes_{t=0}^n \Xi_{t}) \times \Gamma_s)
\rangle {}_{{}_{L^\infty(\Omega_0)}}
\nonumber \\
=
&
{}_{{}_{L^1(\Omega_1)}} \langle   \Phi^{0,1}_* (F_0(\Xi_0) \rho_0), 
{\widehat F}_{1}((\bigtimes_{t=1}^n \Xi_{t}) \times \Gamma_s)
\rangle {}_{{}_{L^\infty(\Omega_1)}}
\label{eq10}
\intertext{
\begin{itemize}
\item[(I)]and,
putting 
$\widetilde{\rho}_0=F_0(\Xi_0) \rho_0$
(or, exactly, its normalization,
i.e.,
$\widetilde{\rho}_0=\lim_{\Xi_0 \to x_0}
\frac{F_0(\Xi_0) \rho_0}{\int_{\Omega_0}{F_0(\Xi_0) \rho_0} d\omega_0}$)
, 
$\widetilde{\rho}_1=F_1(\Xi_1) \Phi^{0,1}_* (\widetilde{\rho}_0)$, 
$\widetilde{\rho}_2=F_2(\Xi_2) \Phi^{1,2}_* ( \widetilde{\rho}_1)$, 
$\cdots$ ,
$\widetilde{\rho}_{s-1}=F_{s-1}(\Xi_{s-1}) \Phi^{s-2,s-1}_* (\widetilde{\rho}_{s-2})$, 
we see that
\end{itemize}}
(\ref{eq10})
=
&
{}_{{}_{L^1(\Omega_1)}} \langle   \Phi^{0,1}_* (\widetilde{\rho}_0), 
{\widehat F}_{1}((\bigtimes_{t=1}^n \Xi_{t}) \times \Gamma_s)
\rangle {}_{{}_{L^\infty(\Omega_1)}}
\nonumber \\
=
&
{}_{{}_{L^1(\Omega_2)}} \langle   \Phi^{1,2}_* ( \widetilde{\rho}_1), 
{\widehat F}_{2}((\bigtimes_{t=2}^n \Xi_{t}) \times \Gamma_s)
\rangle {}_{{}_{L^\infty(\Omega_2)}}
\nonumber \\
&
\cdots \cdots
\nonumber \\
=
&
{}_{{}_{L^1(\Omega_{s+1})}} \langle   
\Phi^{s,s+1}_* ( \widetilde{\rho}_{s}), 
{\widehat F}_{s+1}((\bigtimes_{t=s+1}^n \Xi_{t}) \times \Gamma_s)
\rangle {}_{{}_{L^\infty(\Omega_{s+1})}}
\nonumber \\
=
&
{}_{{}_{L^1(\Omega_{s})}} \langle   
\Phi^{s-1,s}_* ( \widetilde{\rho}_{s-1}), 
{\widehat F}_{s}((\bigtimes_{t=s}^n \Xi_{t}) \times \Gamma_s)
\rangle {}_{{}_{L^\infty(\Omega_{s})}}
\nonumber \\
=
&
{}_{{}_{L^1(\Omega_{s})}} \langle   
\Phi^{s-1,s}_*  (\widetilde{\rho}_{s-1}), 
F_{s} (\Xi_{s})G_{s} (\Gamma_{s})
\Phi^{s,s+1}{\widehat F}_{s+1}(\bigtimes_{t=s+1}^n \Xi_{t})
\rangle {}_{{}_{L^\infty(\Omega_{s})}}
\nonumber \\
=
&
{}_{{}_{L^1(\Omega_{s})}} \langle 
\Big(
F_{s} (\Xi_{s})
\Phi^{s,s+1}{\widehat F}_{s+1}(\bigtimes_{t=s+1}^n \Xi_{t})
\Big)
\Big(  
\Phi^{s-1,s}_* 
(\widetilde{\rho}_{s-1})
\Big), 
G_{s} (\Gamma_{s})
\rangle {}_{{}_{L^\infty(\Omega_{s})}}
\label{eq11}
\end{align}

Thus, we see
\begin{align}
[B_{\widehat{\mathsf{O}}_{0} }^s (\bigtimes_{t \in T} \{x_t \})](\rho_0)
=
\lim_{ \Xi_t \to x_t \;\;(t\in T)}
\frac{
\Big(
F_{s} (\Xi_{s})
\Phi^{s,s+1}{\widehat F}_{s+1}(\bigtimes_{t=s+1}^n \Xi_{t})
\Big)
\times
\Big(  
\Phi^{s-1,s}_*  \widetilde{\rho}_{s-1})
\Big)
}
{\int_{\Omega_0}  [{\widehat F}_{0}(\bigtimes_{t=0}^n \Xi_{t}) ](\omega_0)\;\;
\rho_0 (\omega_0 )   d\omega_0}
\label{eq12}
\end{align}

\subsection{The calculation of
$\rho_s
=
\Phi^{s-1,s}_* (\widetilde{\rho}_{s-1})$
in (\ref{eq12})
}

\par
\noindent
\bf
Lemma 1.
\sl
It holds that
\begin{itemize}
\rm
\item[(J$_1$)]
\sl
$
\int_{- \infty}^{\infty}
\frac{1}{\sqrt{2 \pi} A}
\exp[
-\frac{(x -By )^2}{2A^2}
]
\frac{1}{\sqrt{2 \pi} C}
\exp[
-\frac{(y- D)^2}{2C^2}
]
dy
=
\frac{1}{\sqrt{2 \pi} \sqrt{A^2 + B^2C^2}}
\exp[
-\frac{(x- BD)^2}{2(A^2 + B^2C^2)}
]
$
\rm
\item[(J$_2$)]
\sl
$
\exp[
-\frac{(A\omega -B )^2}{2E^2}
]
\exp[
-\frac{(C \omega- D)^2}{2F^2}
]
\approx
\exp
[
-\frac{1}{2}
(
\frac{A^2 F^2+ C^2 E^2}{E^2 F^2}
)
\Big(\omega -
\frac{(
{AB}{F^2}
+
{CD}{E^2}
)}{({A^2}{F^2} + {C^2}{E^2})}
\Big)^2
]
$
\end{itemize}
where the notation {\lq\lq}$\approx${\rq\rq} means as follows:
$$
"f(\omega) \approx g(\omega )"
\Longleftrightarrow
\mbox{"there exists a positive $K$ such that $f(\omega) = K g(\omega )\;\;
(\forall \omega \in \Omega)$"}
$$
\rm
\par
\noindent
The proof is elementary, and thus, it is omitted.

\par
\vskip1.0cm

\par
\noindent
We see, by (\ref{eq6}) and (I), that
\begin{align}
\widetilde{\rho}_0(\omega_0)=&\lim_{\Xi_0 \to x_0} \frac{F(\Xi_0) \rho_0}{\int_{\mathbb R} F(\Xi_0) \rho_0 d \omega_0}
\nonumber \\
\approx
&
\frac{1}{\sqrt{2 \pi} {q}_0}
\exp[
-\frac{(x_0 -c_0 \omega_0 -d_0 )^2}{2{q}_0^2}
]
\frac{1}{\sqrt{2  \pi} \sigma_0}
\exp[-
\frac{(\omega_0-\mu_0)^2}{2 \sigma_0^2}
]
\nonumber \\
\approx
&
\frac{1}{\sqrt{2  \pi} \widetilde{\sigma}_0}
\exp[-
\frac{(\omega_0-\widetilde{\mu}_0)^2}{2 \widetilde{\sigma}_0^2}
]
\label{eq13}
\end{align}
where
\begin{align}
\widetilde{\sigma}_0^2=\frac{q_0^2 \sigma_0^2}{q_0^2 + c_0^2 \sigma_0^2},
\quad
\widetilde{\mu}_0
=
\mu_0 + \widetilde{\sigma}_0^2 ( 
\frac{c_0}{q_0^2})(x_0-d_0- c_0 \mu_0)
\label{eq14}
\end{align}

\par
\noindent
Further,
the (J$_1$) in Lemma 1 and (\ref{eq9})
imply that
\begin{align}
\rho_1(\omega_1)
&
=
[\Phi_*^{0,1}
\widetilde{\rho}_0](\omega_1) 
\nonumber \\
&
=
\int_{- \infty}^{\infty}
\frac{1}{\sqrt{2 \pi} r_1}
\exp[
-\frac{(\omega_1 - a_1 \omega_{0} -b_1 )^2}{2r_1^2}
]
\frac{1}{\sqrt{2  \pi} \widetilde{\sigma}_0}
\exp[-
\frac{(\omega_0-\widetilde{\mu}_0)^2}{2 \widetilde{\sigma}_0^2}
]
d \omega_0
\nonumber \\
&
=
\frac{1}{\sqrt{2  \pi}{\sigma}_1}
\exp[-
\frac{(\omega_1-{\mu}_1)^2}{2 {\sigma_1}^2}
]
\label{eq15}
\end{align}
where
\begin{align}
{\sigma}_1^2=\ a_1^2 \widetilde{\sigma}_0^2 + r_1^2,
\quad
{\mu}_1
=
a_1 \widetilde{\mu}_0 + b_1 
\label{eq16}
\end{align}

\par
\noindent
Thus,
we see, by (J$_2$) in Lemma 1, that
\begin{align}
\widetilde{\rho}_{t-1}(\omega_{t-1}) 
=
&
\lim_{\Xi_{t-1} \to x_{t-1}} \frac{F(\Xi_{t-1}) \rho_{t-1}}{\int_{\mathbb R} F(\Xi_{t-1}) \rho_{t-1} d \omega_{t-1}}
\nonumber \\
\approx
&
\frac{1}{\sqrt{2 \pi} {q}_{t-1}}
\exp[
-\frac{(x_{t-1} -c_{t-1} \omega_{t-1} - d_{t-1})^2}{2{q}_{t-1}^2}
]
\frac{1}{\sqrt{2  \pi} \sigma_{t-1}}
\exp[-
\frac{(\omega_{t-1}-\mu_{t-1})^2}{2 \sigma_{t-1}^2}
]
\nonumber \\
\approx
&
\frac{1}{\sqrt{2  \pi} \widetilde{\sigma}_{t-1}}
\exp[-
\frac{(\omega_{t-1}-\widetilde{\mu}_{t-1})^2}{2 \widetilde{\sigma}_{t-1}^2}
]
\label{eq17}
\end{align}
where
\begin{align}
\widetilde{\sigma}_{t-1}^2 &=\frac{q_{t-1}^2 \sigma_{t-1}^2}{q_{t-1}^2 + c_{t-1}^2 \sigma_{t-1}^2}
=
\sigma_{t-1}^2
\frac{q_{t-1}^2 + c_{t-1}^2 \sigma_{t-1}^2 +q_{t-1}^2 -q_{t-1}^2 - c_{t-1}^2 \sigma_{t-1}^2}{q_{t-1}^2 + c_{t-1}^2 \sigma_{t-1}^2}
\nonumber \\
&
=
\sigma_{t-1}^2
(1-
\frac{  c_{t-1}^2 \sigma_{t-1}^2}{q_{t-1}^2 + c_{t-1}^2 \sigma_{t-1}^2}
)
\quad
\nonumber \\
\widetilde{\mu}_{t-1}
&
=
\mu_{t-1} + \widetilde{\sigma}_{t-1}^2 ( 
\frac{c_{t-1}}{q_{t-1}^2})(x_{t-1}- c_{t-1} \mu_{t-1})
\label{eq18}
\end{align}

\par
\noindent
Further,
we see, by (J$_1$) in Lemma 1, that
\begin{align}
\rho_t(\omega_{t})
&
=
[\Phi_*^{{t-1},{t}}
\widetilde{\rho}_{t-1}](\omega_{t}) 
\nonumber \\
&
\approx
\int_{- \infty}^{\infty}
\frac{1}{\sqrt{2 \pi} r_{t}}
\exp[
-\frac{(\omega_{t} - a_{t} \omega_{{t-1}} -b_{t} )^2}{2r_{t}^2}
]
\frac{1}{\sqrt{2  \pi} \widetilde{\sigma}_{t-1}}
\exp[-
\frac{(\omega_{t-1}-\widetilde{\mu}_{t-1})^2}{2 \widetilde{\sigma}_{t-1}^2}
]
d \omega_{t-1}
\nonumber \\
&
\approx
\frac{1}{\sqrt{2  \pi} {\sigma}_{t}}
\exp[-
\frac{(\omega_{t}-{\mu}_{t})^2}{2 {\sigma_{t}}^2}
]
\label{eq19}
\end{align}
where
\begin{align}
{\sigma}_t^2=\ a_t^2 \widetilde{\sigma}_{t-1}^2 + r_t^2,
\quad
{\mu}_t
=
a_t \widetilde{\mu}_{t-1} + b_t 
\label{eq20}
\end{align}
%
%
%
%
\par
\noindent
Summing up the above (\ref{eq13})--(\ref{eq20}), we see:
{\small
\begin{align*}
\overset{\mbox{}}{\underset{\mu_0, \sigma_0}{\fbox{\mbox{$\rho_0$}}}}
\xrightarrow[(\ref{eq14})]{\mbox{$x_0$}}
\underset{\widetilde{\mu}_0, \widetilde{\sigma}_0}{\fbox{\mbox{$\widetilde{\rho}_0$}}}
\xrightarrow[(\ref{eq16})]{\Phi_*^{0,1}}
\overset{\mbox{}}{\underset{\mu_1, \sigma_1}{\fbox{\mbox{$\rho_1$}}}}
\xrightarrow[]{\mbox{$x_1$}}
\cdots
\xrightarrow[]{\Phi_*^{t-2,t-1}}
\overset{\mbox{}}{\underset{\mu_{t-1}, \sigma_{t-1}}{\fbox{\mbox{$\rho_{t-1}$}}}}
\xrightarrow[(\ref{eq18})]{\mbox{$x_{t-1}$}}
\underset{\widetilde{\mu}_{t-1}, \widetilde{\sigma}_{t-1}}{\fbox{\mbox{$\widetilde{\rho}_{t-1}$}}}
\xrightarrow[(\ref{eq20})]{\Phi_*^{t-1,t}}
\overset{\mbox{}}{\underset{\mu_t, \sigma_t}{\fbox{\mbox{$\rho_t$}}}}
\xrightarrow[]{\mbox{$x_{t+1}$}}
\cdots
\xrightarrow[]{\Phi_*^{s-1,s}}
\overset{\mbox{}}{\underset{\mu_s, \sigma_s}{\fbox{\mbox{$\rho_s$}}}}
\end{align*}
}
And thus,
we get
\begin{align}
\rho_s
=
\Phi^{s-1,s}_* (\widetilde{\rho}_{s-1})
\label{eq21}
\end{align}
in (\ref{eq12}).
\subsection{
The calculation of
$
\Big(
F_{s} (\Xi_{s})
\Phi^{s,s+1}{\widehat F}_{s+1}(\bigtimes_{t=s+1}^n \Xi_{t})
\Big)
$
in (\ref{eq12})
}

%

\par
Put
\begin{align}
\widetilde{f}_{x_n}(\omega_n)
&
=
\frac{1}{\sqrt{2 \pi} {q}_n}
\exp[
-\frac{(x_n -c_n \omega_n - d_n)^2}{2{q}_n^2}
]
\nonumber \\
&
\approx
\exp[
-\frac{(c_n \omega_n -( x_n-d_n))^2}{2{q}_n^2}
]
\equiv
\exp[
-
\frac
{1}{2}
\Big(\widetilde{u}_n \omega_n -\widetilde{v}_n \Big)^2
]
\label{eq22}
\end{align}
where it is assumed that 
$c_n$, $d_n$ and $q_n$ are known
$(t \in T)$. And thus, put
\begin{align}
{\widetilde{u}_n}
=\frac{c_n}{q_n}
,\quad
{\widetilde{v}_n}
=
\frac{x_n - d_n}{q_n}
\label{eq23}
\end{align}

\par
\noindent
And further, Lemma 1 implies that the causal operator
$\Phi^{t-1.t}:L^\infty (\Omega_t) \to L^\infty ( \Omega_{t-1})$
is defined by
\begin{align}
&
f_{t-1}(\omega_{t-1})
=
[\Phi^{t-1,t} \widetilde{f}_{x_t}](\omega_{t-1} )
\nonumber \\
\approx
&
\int_{- \infty}^{\infty}
\frac{1}{\sqrt{2 \pi} r_t}
\exp[
-\frac{(\omega_t - a_t \omega_{t-1} -b_t )^2}{2r_t^2}
]
\exp[
-\frac{({\widetilde{u}}_t \omega_t-{\widetilde{v}_t})^2}{2}
]
d \omega_{t}
\nonumber \\
\approx
&
\exp[
-\frac{1}{2}
\Big(
\frac{\widetilde{v}_t}{\sqrt{1 +r_t^2 {\widetilde{u}}_t^2 }}
-
\frac{{\widetilde{u}}_t(a_t \omega_{t-1} + b_t)}{\sqrt{1 +r_t^2 {\widetilde{u}}_t^2}}
\Big)^2
]
\approx
\exp[
-\frac{1}{2}
\Big({{u}}_{t-1} \omega_{t-1}-{{v}_{t-1}}
\Big)^2
]
\label{eq24}
\end{align}
where
\begin{align}
&
u_{t-1}=-\frac{a_t {\widetilde u}_t}{\sqrt{1 +r_t^2 {\widetilde{u}}_t^2 }},
\quad
v_{t-1}=\frac{  b_t {\widetilde u}_t - {\widetilde v}_t }{\sqrt{1 +r_t^2 {\widetilde{u}}_t^2 }}
\label{eq25}
\end{align}

\par
\noindent
And also, Lemma 1 implies that
\begin{align}
\widetilde{f}_{x_{t-1}}(\omega_{t-1})
&
=
\exp[
-\frac{(c_{t-1} \omega_{t-1} + d_{t-1}-x_{t-1})^2}{2{q}_{t-1}^2}
]
\exp[
-\frac{({{u}}_{t-1} \omega_{t-1}-{{v}_{t-1}})^2}{2}
]
\nonumber \\
&
\approx
\exp
[
-\frac{1}{2}
(
\frac{c_{t-1}^2 + u_{t-1}^2 q_{t-1}^2}{q_{t-1}^2 }
)
\Big(\omega_{t-1} -
\frac{
{c_{t-1}(d_{t-1}-t_{t-1})}
+
{u_{t-1}v_{t-1}}{q_{t-1}^2}
}{{c_{t-1}^2} + {u_{t-1}^2}{q_{t-1}^2}}
\Big)^2
]
\nonumber \\
&
\approx
\exp[
-\frac{1}{2}
\Big({\widetilde{u}}_{t-1} \omega_{t-1}-{\widetilde{v}_{t-1}}
\Big)^2
]
\label{eq26}
\end{align}
where
\begin{align}
{{\widetilde{u}}_{t-1}}
=
\frac{\sqrt{c_{t-1}^2 + u_{t-1}^2 q_{t-1}^2}}{q_{t-1} },
\;\;
{{\widetilde{v}}_{t-1}}
=
\frac{
{c_{t-1}(d_{t-1}-t_{t-1})}
+
{u_{t-1}v_{t-1}}{q_{t-1}^2}
}{q_{t-1} \sqrt{{c_{t-1}^2} + {u_{t-1}^2}{q_{t-1}^2}}}
\label{eq27}
\end{align}

Summing up the above (\ref{eq22})-(\ref{eq27}), we see:
{\scriptsize
\begin{align*}
{
\overset{\widetilde{u}_s, \widetilde{v}_s}{
\underset{\widetilde{w}_s}{\fbox{\mbox{$\widetilde{f}_{x_s}$}}}
}}
\xleftarrow[]{\mbox{$x_s$}}
\cdots
\xleftarrow[]{\Phi^{t-2,t-1}}
{
\overset{\widetilde{u}_{t-1}, \widetilde{v}_{t-1}}{
\underset{\widetilde{w}_{t-1}}{
{{\fbox{\mbox{$\widetilde{f}_{x_{t-1}}$}}}}}}}
\xleftarrow[(\ref{eq27})]{{x_{t-1}}}
{
\overset{{u}_{t-1}, {v}_{t-1}}{
\underset{w_{t-1}}{
\fbox{\mbox{${f}_{t-1}$}}}}
}
\xleftarrow[(\ref{eq25})]{\Phi^{t-1,t}}
{
\overset{\widetilde{u}_t, \widetilde{v}_t}{
\underset{{\widetilde{w}_t}}{
\fbox{\mbox{$\widetilde{f}_{x_t}$}}}}
}
\xleftarrow[]{{x_{t}}}
\cdots
\xleftarrow[]{{x_{n-1}}}
{
\overset{{u}_{n-1}, {v}_{n-1}}{
\underset{{w}_{n-1}}{
\fbox{\mbox{${f}_{n-1}$}}}}
}
\xleftarrow[]{\Phi^{n-1,n}}
{
\overset{{\widetilde{u}_n \widetilde{v}_n}}{\underset{{\widetilde{w}_n}}{\fbox{\mbox{$\widetilde{f}_{x_n}$=(\ref{eq23})}}}}
}
\end{align*}
}
And thus,
we get
\begin{align}
\widetilde{f}_{x_s}
\approx
\lim_{\Xi_t \to x_t \;( t \in \{s.s+1, \cdots, n\})}
\frac{
\Big(
F_{s} (\Xi_{s})
\Phi^{s,s+1}{\widehat F}_{s+1}(\bigtimes_{t=s+1}^n \Xi_{t})
\Big)}
{
\|
F_{s} (\Xi_{s})
\Phi^{s,s+1}{\widehat F}_{s+1}(\bigtimes_{t=s+1}^n \Xi_{t})
\Big)
\|_{L^\infty (\Omega_s )}
}
\label{eq28}
\end{align}
in (\ref{eq12})
\par
\noindent
\par
\par
After all,
we solve the problem (H),
that is,
\begin{itemize}
\item[(H$'$)]
Assume that
a measured value
$(x_0, x_2, \cdots, x_n )$
$(\in \bigtimes_{t=0}^n X_t)$
is obtained by the measurement
${\mathsf M}_{L^\infty (\Omega_0)}$
$(\widehat{\mathsf O}_{t_0},$
$
S_{[\ast]}(\rho_0 )
)$.
Let $s(\in T)$ be fixed. Then,
we get
the Bayes-Kalman operator
$[B_{\widehat{\mathsf{O}}_{t_0} }^s(\bigtimes_{t \in T} \{x_t\})]
(\rho_0)$, that is,
$$
\Big([B_{\widehat{\mathsf{O}}_{t_0} }^s(\bigtimes_{t \in T} \{x_t\})]
\rho_0
\Big)(\omega_s)
=
\frac{\widetilde{f}_{x_s}(\omega_s ) \cdot \rho_s (\omega_s )}{
\int_{-\infty}^{\infty}
\widetilde{f}_{x_s}(\omega_s ) \cdot \rho_s (\omega_s )
d \omega_s
}
\quad
(
\forall \omega_s \in \Omega_s 
)
$$
where
$\rho_s$
in (\ref{eq21})
and
$\widetilde{f}_{x_s}$
in
(\ref{eq28})
can be iteratively calculated as mentioned in this section.
\end{itemize}
\par
\noindent
\vskip0.5cm
\par
\noindent
{\it Remark 2.}
The following clssification is usual
\begin{itemize}
\item[(K$_1$)]
Smoothing: in the case that $0 \le s < n$
\item[(K$_2$)]
Filter:
in the case that $s= n$
\item[(K$_3$)]
Prediction:
in the case that $s= n$ and,
for any $m$
such that $n_0 \le m < n$,
the observable $(X_m, {\mathcal F}_m, F_m )=
(\{1\}, \{\emptyset ,\{1 \} \}, F_m )$
is defined by
$F_m(\emptyset )\equiv 0$,
$F_m(\{ 1 \} )\equiv 1$,
\end{itemize}

\section{Conclusions}
\par
\noindent
\par
We conclude as follows.
\begin{itemize}
\item[(L)]
Mathematically, quantum language may be somewhat more difficult than usual statistics. In fact, quantum language can not be understood without the knowledge of
functional analysis.
However,
we think that the calculation of Kamlan filter is more understandable in terms of quantum language than in terms of usual statistics.
Therefore we believe that
quantum language is future statistics,
i.e.,
statistics will develop into quantum language. 
\end{itemize}

\par
We hope that our proposal will be discussed and examined from various view-points.

\rm
\par


\rm
\par
\renewcommand{\refname}{
\large 
References}
{
\small

\normalsize
}

\end{document}